\documentclass[11pt]{article}
\usepackage{graphicx}
\usepackage{amssymb, amsmath, latexsym}
\usepackage{epstopdf}
\usepackage{amscd}
\newtheorem{theorem}{Theorem}
\newtheorem{proposition}{Proposition}
\newtheorem{corollary}[theorem]{Corollary}
\newtheorem{definition}{Definition}
\newtheorem{lemma}{Lemma}
\newtheorem{prop}{Proposition}

\newcommand{\spine}[2]{{\mathfrak C}_{#1,#2}}
\newcommand{\dspine}[2]{{\mathfrak C}^{#1,#2}}

\newcommand{\proof}[1]{\textbf{Proof: }{#1} $\Box$}

\newcommand{\Q}{\mathbb Q}

\newcommand{\Out}[1]{\operatorname{Out}(F_{#1})}
\newcommand{\Aut}[1]{\operatorname{Aut}(F_{#1})}
\newcommand{\bdry}{\partial}

\title{Morita classes in the homology of $\Aut{n}$ vanish after one stabilization
}
\author{James Conant\\ Karen Vogtmann}
\begin{document}
\maketitle
\begin{abstract} 
There is a series of cycles in the rational homology of the groups $\Out{n}$, first discovered by S. Morita, which have an elementary description in terms of finite graphs.  The first two of these give nontrivial homology classes, and it is conjectured that they are all nontrivial. These cycles have natural lifts to the homology of $\Aut{n}$, which is stably trivial by a recent result of Galatius.  We show that in fact a single application of the stabilization map $\Aut{n} \to \Aut{n+1}$ kills the Morita classes, so that they  disappear immediately after they appear.
   
\end{abstract}

\section{Introduction}

In \cite{Mor99} S. Morita   constructed a series of cycles in the rational homology of $\Out{n}$. These homology classes are related to Morita's trace map, defined in \cite{Mor93}, which has interesting connections to a number of different problems in topology.  These include understanding the group of homology cobordism classes of homology cylinders and understanding the arithmetic mapping class group, an extension of the mapping class group by the  Galois group of $\overline \Q/\Q$ (for more information, see Morita's survey article, \cite{Mor05}).  Morita showed that the first of his cycles gives a non-trivial homology class in $H_4(\Out{4};\Q)$, and Conant and Vogtmann \cite{ConVog04}  proved that the second cycle also gives a non-trivial homology class, this time in $H_8(\Out{6};\Q)$.  It is conjectured that in fact all of Morita's cycles are nontrivial in homology, a  conjecture which would have intriguing implications for the related problems.  

The homology $H_4(\Out{4};\Q)\cong H_4(\Aut{4};\Q)$ was known to be 1-dimensional \cite{HatVog98:2}, so that Morita's class generates the entire group.  Recently Ohashi calculated that $H_8(\Out{6};\Q)$ is also 1-dimensional \cite{Oha06}, so again the Morita cycle generates the entire homology. 
This  leads to the more speculative conjecture that Morita's cycles, together with cycles obtained by a natural generalization of Morita's construction, in fact generate all of the rational homology of $\Out{n}$.  This homology has been completely computed up to dimension 7, but $H_8$ is not fully known. For example, $H_8(\Out{7};\Q)$ is not known, and indeed two Morita cycles give potential classes there.

Morita's cycles have natural lifts to cycles for $\Aut{n}$, which we also call Morita cycles.  We note that Gerlits \cite{Ger} found a rational class in $H_7(\Aut{5};\Q)$ that, for degree reasons, cannot be a Morita cycle, so that  the rational homology of $\Aut{n}$ is not given entirely by these classes.  

 Hatcher and Vogtmann \cite{HatVog98:3} showed that the natural map $\Aut{n}\to \Aut{n+1}$ induces an isomorphism $H_i(\Aut{n};\Q)\to H_i(\Aut{n+1};\Q))$  for $n\geq 5i/4$.  In addition, the map $H_i(\Aut{n})\to H_i(\Out{n})$ is an isomorphism for $n>2i+4$ (\cite{HatVog04,HatVogWah06}), so that the homology of $\Out{n}$ is also independent of $n$ for $n$ sufficiently large.  The purpose of this article is to show that the  Morita cycles  in the homology of $\Aut{n}$  vanish immediately after they appear:  the   $i$-th Morita class is a cycle 
in $Z_{4k}(\Aut{2k+2})$, and we show that the image of this cycle in $Z_{4k}(\Aut{2k+3})$ is a boundary, and hence the image in $Z_{4k}(\Out{2k+3})$ is also a boundary.   

In a recent paper, S. Galatius constructed an infinite loop space  which computes the stable homology of $\Aut{n}$ (and $\Out{n}$)  and then proved that this space is homotopy equivalent to $\Omega^\infty S^\infty$ \cite{Gal06}.  In particular, this shows that the  stable rational homology   of $\Aut{n}$ and $\Out{n}$ is zero.  Thus the Morita classes, which are in the unstable range, must eventually vanish under the stabilization maps $\Aut{n}\to \Aut{n+1}\to \Aut{n+2}\to\cdots$.  Our theorem shows that this happens after a single stabilization.  If the Morita classes do in fact generate the homology, then this homology is ephemeral indeed, and the homology stabilizes earlier than stated in the known stability theorem.

{\bf Acknowledgments:} The first author is supported by NSF grant DMS 0604351. The second author is supported by NSF grant DMS 0204185.

\section{Chain and cochain complexes}

In order to compute rational homology of $\Aut{n}$ and $\Out{n}$ we need to produce  suitable chain complexes.  In this section we describe how this is done.

The spine  of Outer space is a locally finite, contractible simplicial complex $K_n$ on which the group $\Out{n}$ acts cocompactly with finite stabilizers (see, e.g., \cite{Vog02}).   Thus the rational homology of the quotient $K_n/\Out{n}$ is equal to the rational homology of $\Out{n}$.  Vertices of $K_n$ are ``minimal marked graphs"  $(g,G)$.   Here the {\it marking} $g$ is a homotopy equivalence from a fixed standard rose $R_n$ to the graph $G$, and a graph is {\it minimal} if it is connected with no univalent or bivalent vertices. If the graphs have specified basepoints, the analogous spine $L_n$ has an action of $\Aut{n}$.  

In \cite{HatVog98:2} it was shown that   
the spine can be given the structure of a cube complex, and a chain complex was constructed to compute the homology of the quotient space. For the convenience of the reader, we briefly recall this construction.  The spine has  one cube $(g, G ,F)$ for every marked graph $(g,G)$  and forest $F\subset G$, where a forest is a subgraph with no cycles.     The dimension of the cube is equal to the number of edges in $F$, and an orientation on the cube is specified by ordering the edges of $F$, up to even permutation.   The codimension 1 faces of $(G,F)$ are of two types, obtained from $(G,F)$ by either removing an edge from $F$ or by collapsing an edge of $F$; in either case the ordering on the remaining edges of $F$ determines an orientation of the face.       

The action of $\Out{n}$ on $K_n$ (or $\Aut{n}$ on $L_n$) simply changes the marking $g$.  This action preserves the cube-complex structure, and is transitive when restricted to a given isomorphism class of pairs $(G,F)$. 
 The quotient $Q_n$ by this action is  a cell complex, with one cell for each isomorphism class of pairs $(G,F)$.  It is no longer a cube complex because faces of a cube may be folded or identified with other faces under the action. The rational homology of $Q_n$ is computed by the chain complex $\mathfrak C_*$, which is the quotient of the vector space spanned by isomorphism classes of oriented pairs $(G,F)$ modulo the relations $(G, F)+(G,-F)=0$.  
  The boundary map is given by the formula
$$\bdry(G,F)=\sum_{e_i\in F} (-1)^i \big((G, F-e_i) - (G_{e_i},F_{e_i})\big),$$
where $(G_{e_i},F_{e_i})$ denotes the result of contracting the edge $e_i$.

The $k$-chains $\mathfrak C_k$ can be decomposed further according to the number of vertices of $G$, i.e. 
$$\mathfrak C_k=\bigoplus_{p+q=k} \mathfrak C_{p,q}$$
where $\mathfrak C_{p,q}$ has one generator for each pair $(G,F)$ such that $G$ has $q$ vertices.  Since $G$ has no univalent or bivalent vertices and $\chi(G)=1-n$, we have $1\leq q\leq 2n-2$ and $2-2n\leq p\leq -1$.  The boundary operator $d\colon \mathfrak C_k\to\mathfrak C_{k-1}$ is the sum of horizontal boundary operators $$d_R\colon \mathfrak C_{p,q}\to \mathfrak C_{p-1,q}$$
corresponding to faces obtained by removing an edge from $F$ and vertical operators
$$d_C\colon \mathfrak C_{p,q}\to \mathfrak C_{p,q-1}$$
corresponding to faces obtained by collapsing an edge of $F$.  Thus $\{\mathfrak C_{p,q}, d_R, d_C\}$ is a second-quadrant double complex so that the horizontal and vertical filtrations of $\{\mathfrak C_{p,q}, d_R, d_C\}$ give rise to second-quadrant spectral sequences converging to the rational homology of $Q_n$.

There is also a dual complex $\displaystyle \mathfrak C^*=\oplus_{p,q}\mathfrak C^{p,q}$  which computes  the cohomology of $Q_n$.  This arises in the same way, except that the coboundary operator is defined on a cube as the sum of cubes for which the original cube is a codimension $1$ face. As before, the coboundary operator decomposes as a sum of
the coboundary operators $\delta_R$ and $\delta_C$, which add an edge to the forest in all ways and expand a vertex into an edge in all ways, respectively.  If $|F|$ denotes the number of edges in a forest $F$, then 
$$\delta_R(G,F)=(-1)^{|F|+1}\sum_e (G, F\cup e)$$
where $e$ runs over all edges $e\in G\setminus F$ with $F\cup e$ a forest and
$$\delta_C(G,F)=(-1)^{|F|}\sum_P (G^P,F^P)$$
where $P$ runs over all partitions of the incoming edges at each vertex into two sets, each of cardinality at least $2$.  Each such partition can be used to blow up the vertex into an edge $e_P$ of the forest, and $(G^P, F^P)$ denotes the resulting forested graph.  The orientations on $(G,F\cup e)$ and $(G^P,F^P)$ are determined by placing the new edges $e$ and $e_P$   last in the ordering of the forest edges.

{\bf Remark:} The signs in the above formulas arise in the following way. Every edge  $e$ of the forest corresponds to an axis of the cube, which starts at the face where the edge is collapsed, and ends at the opposite face where the edge is removed from the forest. This determines a unit tangent vector $\vec{v}_e$ in the tangent bundle of the cube, the fibres of which are all canonically identified. Thus an ordering of the edges gives rise to an orientation $or=\vec{v}_{e_1}\wedge\cdots\wedge \vec{v}_{e_{|F|}}$. In (co)homology calculations, the induced orientation when passing to a face is given by rewriting the orientation as $or =\vec{n}\wedge \nu$, where $\vec{n}$ is the unit inward-pointing normal to the face. One then takes $\nu$ as the induced orientation. The signs in $d_C,d_R,\delta_C,\delta_R$ can all be easily calculated using this formalism. For example, the sign for $\delta_R$ is computed by analyzing the sign from passing from the cube $(G,F\cup e)$ to the face $(G,F)$. Let $\mathbf{v_F}$ be the orientation of this face given by $\mathbf{v_F}=\vec{v}_{e_1}\wedge\cdots\wedge \vec{v}_{e_{|F|}}$ coming from the edge ordering. Then the orientation of the whole cube is given by $\mathbf{v_{F}}\wedge \vec{v}_e$. The inward pointing normal to the face is $-\vec{v}_e$, and so $\mathbf{v_{F}}\wedge \vec{v}_e=(-1)^{|F|+1}(-\vec{v}_e)\wedge \mathbf{v_F}$, giving rise to the sign $(-1)^{|F|+1}$.

The  complexes  $\mathfrak C_*$ and $\mathfrak C^*$ are related by a bilinear form $$\langle\cdot,\cdot\rangle\colon \mathfrak C^*\otimes\mathfrak C_*\to \mathbb Q,$$ 
defined on the level of oriented pairs $X=(G,F)$ as follows: $\langle X,Y\rangle$ is equal to zero unless $X$ and $Y$ are isomorphic as pairs  and neither has an orientation reversing automorphism.  In this case  \[
\langle X,Y\rangle=
\begin{cases} 
|\operatorname{Iso}(X,Y)| & \text{if $X\cong Y$,}\\
-|\operatorname{Iso}(X,-Y)| & \text{if $X\cong -Y$}
\end{cases}
\]
where $\operatorname{Iso}(X,Y)$ is the number of (orientation-preserving) isomorphisms between $X$ and $Y$.

\begin{proposition}
$\langle (\delta_R+\delta_C)X,Y\rangle =\langle X,(d_R+d_C)Y\rangle$.
\end{proposition}

\proof{
We actually prove that both $\langle \delta_R X,Y\rangle = \langle X,d_RY\rangle$ and 
$\langle \delta_C X,Y\rangle = \langle X,d_CY\rangle$ for all oriented pairs $X=(G,F)$ and $Y=(G',F')$.  

Suppose $X$ has an orientation reversing automorphism. Then $X$ and hence $\delta_RX$ and $\delta_CX$ are all zero in $\mathfrak C^*$. Thus the two equations trivially hold. Similarly, the equations trivially hold if $Y$ has orientation reversing automorphisms. So we may assume that neither $X$ nor $Y$ has an orientation-reversing automorphism.

We first establish $\langle \delta_R X,Y\rangle = \langle X,d_RY\rangle$.
Consider the sets
\begin{align*}
A_+&=\{ (e,\phi) | e\in G\setminus F, \phi\colon (G,F\cup e)\overset{\cong}{\to}(G',F')\}\\
A_-&=\{ (e,\phi) | e\in G\setminus F, \phi\colon (G,F\cup e)\overset{\cong}{\to}-(G',F')\}
\end{align*}
where $\phi$ ranges over all orientation-preserving isomorphisms.
Then $\langle \delta_R X,Y\rangle=(-1)^{|F|+1}(|A_+|-|A_-|)$, for each term of $\delta_R X$ corresponds to a choice of $e$, and the inner product then decides if the resulting graph is isomorphic to $\pm Y$, multiplying by $|\operatorname{Aut}(Y)|$ if it is. This is the same as counting all possible isomorphisms with $\pm Y$, (which all must have the same sign since $Y$ has no orientation-reversing symmetries.)

Similarly define
\begin{align*}
B_+&=\{ (e_i',\psi) | e_i'\in F', \psi\colon (-1)^{|F|+1-i}(G',F'- e_i')\overset{\cong}{\to}(G,F)\}\\
B_-&=\{ (e_i',\psi) | e_i'\in F', \psi\colon (-1)^{|F|+1-i}(G',F'- e_i') \overset{\cong}{\to}-(G,F)\}
\end{align*}
where $\phi'$ ranges over orientation-preserving isomorphisms. Then  $\langle X,d_RY\rangle=(-1)^{|F|+1}(|B_+|-|B_-|)$.
 
We claim that there is a bijection from $A_+$ to $B_+$. Given $(e,\phi)\in A_+$, set  $e_i'=\phi(e)$ and let $\psi$ be equal to the inverse of the restriction of $\phi$ to $(G,F)$.  Since $\phi$ is an isomorphism, the ordering on the edges of $F'$ induced by $\phi$ is compatible with the orientation of $F'$.  In this ordering, $\phi(e)$ is last, i.e.  $i=|F|+1$, so $(-1)^{|F|+1-i}(G',F'- e_i')=(G',F'- e_i')$ and the formula  is correct.  For the inverse map, given $(e_i',\psi)\in B_+$, set $e=\psi^{-1}(e_i')$ and let $\phi$ be equal to the extension of $\psi^{-1}$ to $(G',F')$.  The sign $(-1)^{|F|+1-i}$ comes from the fact that we need to move $e_i'$ to the end of the ordering since the ordering of $F\cup e$ has $e$ last.   
An identical argument gives a bijection from $A_-$ to $B_-$, and we conclude that $\langle \delta_R X,Y\rangle = \langle X,d_RY\rangle$.

 Next we show that $\langle \delta_C X,Y\rangle = \langle X,d_CY\rangle$.  As above, we define sets
 $$A_\pm=\{(P,\phi)| P\text{ is a partition of edges at a vertex}, \phi\colon(G^P,F^P) \overset{\cong}{\to} \pm(G',F')\}$$
 where $\phi$ ranges over all orientation-preserving isomorphisms. Then 
 $\langle \delta_C X,Y\rangle=(-1)^{|F|}(|A_+|-|A_-|).$ Also as above set
 $$B_{\pm}=\{(e_i',\psi) | e_i'\in F', \psi\colon (-1)^{|F|+1-i}(G'_{e_i'},F'_{e_i'}) \overset{\cong}{\to} \pm(G,F)\}$$
 where $\phi'$ ranges over orientation preserving isomorphisms. We still have
 $\langle X,d_CY\rangle=(-1)^{|F|}(|B_+|-|B_-|)$. (Recall that $d_C$ has an extra minus sign.)
 
 As before, we will show that $A_\pm$ has the same size as $B_\pm$. Given $(P,\phi)\in A_\pm$ let $e_i'=\phi(e_P)$ and define $\psi$ as follows. First note that $\phi^{-1}$ gives an isomorphism from $(G',F')$ to $(G^P,F^P)$. Thus we get an induced isomorphism $(G'_{e_i'},F'_{e_i'})\to ((G^P)_{e_P},(F^P)_{e_P})=(G,F)$, which we define to be $\psi$. This gives a map $A_\pm\to B_\pm$. The inverse map is defined analogously.}


\section{Morita Cycles}

In this section we fix an odd-valent graph $\gamma$ together with certain orientation data, and use it  to construct a cycle  $z(\gamma)$ in the chain complex $\mathfrak C_*$.  In the next section we will show that this cycle corresponds to the cocycle $\mu_\gamma$ defined in \cite{ConVog04}.

The idea is to form $z(\gamma)$ by replacing each vertex $x$ of $\gamma$ by an $|x|$-gon $C_x$ (where $|x|$ is the valence of $x$),  i.e. removing $x$ and attaching the edges that used to terminate at $x$ to the vertices of $C_x$.  The isomorphism type of the resulting graph depends on exactly how we reattach the edges to $C_x$, so we make this operation canonical by summing over all possibilities.  With a careful choice of signs, this sum is a cycle; in fact, it is a cycle for both components $d_R$ and $d_C$ of the boundary operator $d$ of $\mathfrak C_*$.   

Here are the details of the construction.

\begin{definition} A graph $\gamma$ is \emph{admissible} if  $\gamma$ is connected and all vertices have odd valence greater than or equal to 3.  
\end{definition}

\begin{definition}
An \emph{orientation at a vertex}  of a graph is an ordering of the half-edges incident to that vertex, up to even permutation. A graph is \emph{vertex-oriented} if every vertex has an orientation
\end{definition}

Let $\gamma$ be an admissible, vertex-oriented graph   For each vertex $x$ of $\gamma$, choose an ordering of the edges at $x$ compatible with the vertex-orientation,  let $C_x$ be an $|x|$-gon, with vertices $v_1,\ldots, v_{|x|}$ ordered cyclically and edges $e_1,\ldots,e_{|x|}$ also ordered cyclically, $e_i$ ending at $v_i$.  

Given any permutation $\sigma_x$ of $\{1,\ldots,|x|\}$ we can form a new graph by removing $x$ then re-attaching the $j$-th edge in the ordering to $v_{\sigma_x(j)}$, for $j=1,\ldots,|x|$.  
Given a collection $\sigma=\{\sigma_x\}$ of permutations, one for each vertex $x$ of $\gamma$,  we can perform this operation at each vertex of $\gamma$ to obtain a new graph $G^\sigma$. 

We will be interested in forests $F$  which are maximal forests  of the subgraph of $G^\sigma$  consisting of the union of the cycles $C_x$.  The intersection $F\cap C_x$ is obtained by deleting a single edge of $C_x$.  If this edge is $e_j$, set $\epsilon(F,x)=(-1)^j$.  

The cycle $z(\gamma)$ is now defined by summing over all collections of permutations $\sigma=\{\sigma_x\}$ and all maximal forests $F$ in $\cup C_x$:
$$z(\gamma)=\sum_\sigma\sum_F\epsilon_\sigma\epsilon_F(G^\sigma,F),$$
where
 $\epsilon_F=\prod_x \epsilon(F,x)$ and $\epsilon_\sigma=\prod_x sign(\sigma_x)$, 

Note that this sum is independent of the choice of representatives for the vertex orientations, since an even permutation of the edges at $x$ does not change the sign of $\sigma_x$.  We also remark that if a vertex $x$ had even valence, then there would be an automorphism of $C_x$ inducing an odd permutation of its vertices, and  the terms in the sum would cancel in pairs, giving $z(\gamma)=0$.  This is the reason for the restriction to odd-valent graphs.

Recall that the differential $d$ in $\mathfrak C_*$ is the sum $d_R+d_C$; we actually show that each of these is zero on $z(\gamma)$.  

\subsection{$d_R(z(\gamma))=0$} 
The ambient $G$ is constant in this computation, so we can suppress it in the notation, and write 
$d_R(F)$ instead of $d_R(G,F)$; thus  $\displaystyle d_R(F)=\sum_{e_i\in F} (-1)^i(F-e_i)$.    
   
In fact, the operator $d_R$ makes sense in a more general context.   An {\it ordered graph}    is a graph $A$ together with an ordering on the edges of $A$.  An ordering on $A$ induces a natural ordering on any subgraph $B$ of $A$.   We say two orderings on $A$ are equivalent if they differ by an even permutation, and let $\mathcal {E}$ denote the vector space spanned by equivalence classes of ordered graphs modulo the relation $(A, ord)=-(A,-ord)$.    An {\it edge-oriented graph} is an equivalence class of ordered graphs. 

If $A$ and $B$ are edge-oriented graphs, choose orderings representing the edge-orientations and define $A\cdot B$ to be the (equivalence class of the) disjoint union of $A$ and $B$, with ordering given by increasing the labels on edges of $B$ by the number $e(A)$  of edges of $A$.  We see immediately that 

\begin{lemma}\label{dot}  
$A\cdot B= (-1)^{e(A)e(B)} B\cdot A$  \end{lemma}

Now define $d_R\colon \mathcal {E}\to \mathcal {E}$ by the formula
$$d_R(A)=\sum_{e_i\in A} (-1)^i (A-e_i)$$
where $e_i$ is the $i$-th edge of $A$.

\begin{lemma}  $d_R^2=0$
\end{lemma}
\proof{  This is the standard argument:
$$d_R^2(A)=\sum_{j<i}(-1)^i(-1)^j (A-e_i-e_j) +\sum_{j>i} (-1)^i(-1)^{j-1}(A-e_i-e_j)$$
so that the terms cancel in pairs.
}

\begin{lemma}   $d_R(A\cdot B)=d_R(A)\cdot B + (-1)^{e(A)} A\cdot d_R(B)$
\end{lemma}
 \proof{This is immediate.}
  
\begin{prop}\label{dR}  Let $G$ be a graph, and $C$ a subgraph which is the disjoint union of polygons $C_1,\ldots, C_k$.  Order the edges of each $C_i$, and let $C=C_1\cdot\ldots\cdot C_k$.   Set  $\epsilon_F = \prod_i(-1)^ {i_j}$, where  if $F$ is obtained by deleting the $i_j$-th edge from $C_i$. Let
 $$Z=\sum_F\epsilon_F(G, F)$$
where  the sum is over all maximal forests $F$ of $C$ and the ordering on $F$ is induced from the ordering on $C=C_1\cdot\ldots\cdot C_k$.  Then    $d_R(Z)=0$.
\end{prop}

\proof{   The terms of 
$d_R(C_1)\cdot d_R(C_2)\cdot \ldots\cdot d_R(C_k)$ are precisely the maximal forests in $C$, with sign $\epsilon_F$.  Thus 
$$d_R(\sum_F \epsilon_F(G, F))=
 d_R(d_R(C_1)\cdot \ldots\cdot d_R(C_k))=\sum_i \pm(d_R(C_1)\cdot \ldots\cdot d_R^2(C_i)\cdot \ldots\cdot d_R(C_k))=0.$$
}

\begin{figure}
\begin{center}
\includegraphics[width=2in]{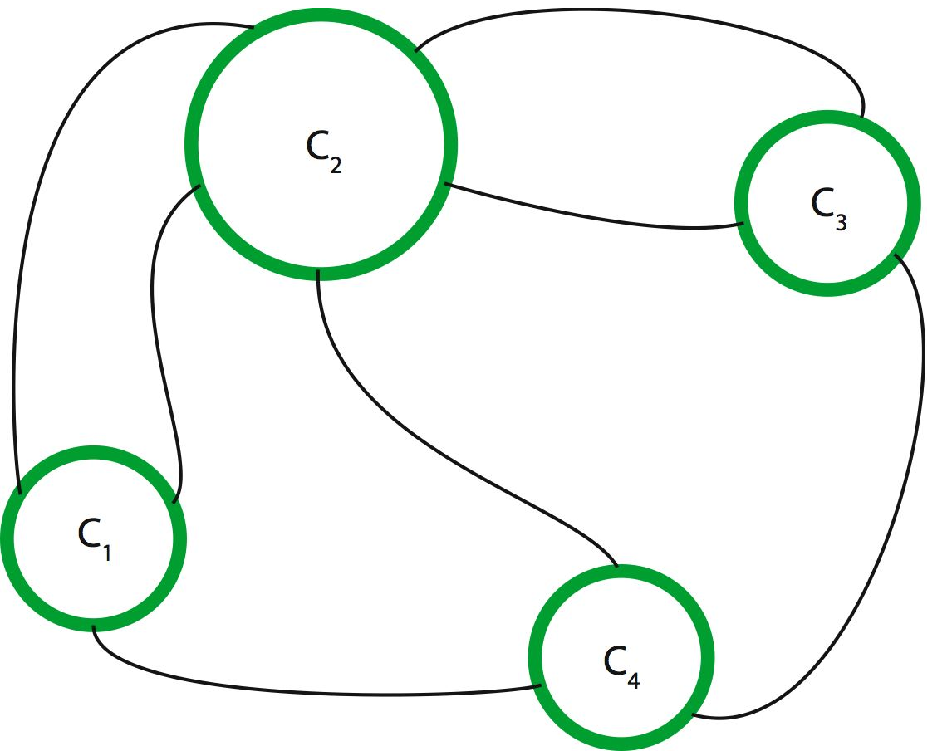}
\caption{The $d_R$-cycle $Z$}\label{zgamma}
\end{center}
\end{figure}

The cycle $Z$ is illustrated in Figure \ref{zgamma}; the forests $F$ are all maximal forests in the green subgraph.  

\begin{corollary}  $d_R(z(\gamma))=0$
\end{corollary}

\proof{ In our construction of $G^\sigma$ we fixed a cyclic ordering of the edges of each $C_{x}$.  Let $Z_\sigma$ denote the corresponding $d_R$-cycle from Proposition~\ref{dR}.  Then 
 $$z(\gamma)=\sum_\sigma\sum_F\epsilon_\sigma\epsilon_F(G^\sigma,F)=\sum_\sigma\epsilon_\sigma Z_\sigma,$$
 so 
$$d_R(z(\gamma))=\sum_\sigma \epsilon_\sigma d_r(Z_\sigma)=0$$
}

\subsection{$d_C(z(\gamma))=0$}  Here, too, we introduce a slightly more general context for future convenience:  

\begin{prop}\label{dc} Let $G$ be a trivalent graph, and $C$ a subgraph which is the disjoint union of polygons $C_1,\ldots, C_k$.  Fix a maximal forest $F$ of $C$, order its edges and set  
$$Z=\sum_\sigma \epsilon_\sigma(G^\sigma, F)$$
where $\epsilon_\sigma$ is the sign of $\sigma$ and $G^\sigma$ is obtained from $G$ by permuting the edges coming into each $C_i$ by the permutation $\sigma_i$. Then   $d_C(Z)=0$.
\end{prop}

\proof{ For this, we simply note that $(G^\sigma_e, F_e)= (G^{\tau}_e, F_e)$ if $\sigma$ and $\tau$ differ only by 
transposing the two edges coming into the vertices of $e$.  Since in this case $\epsilon_\sigma=-\epsilon_\tau$,
all terms of $d_C(z(\gamma))$ cancel in pairs.
}

\begin{corollary}
$d_C(z(\gamma))=0$.
\end{corollary}
\proof{
We have $$d_C(z(\gamma))=d_C(\sum_\sigma\sum_F \epsilon_\sigma \epsilon_F (G^\sigma,F))=\sum_F \epsilon_F\, d_C(\sum_\sigma \epsilon_\sigma(G^\sigma,F))=\sum_F0=0.$$
}

\noindent{\bf Remark.}  
The same construction as above gives a cycle $z(\gamma)$ associated to any graph $\gamma$ with orientations at any (non-empty) subset of its odd-valent vertices. These cycles correspond to the cocycles defined in \cite{ConVog04} which take values in the space of ``AB-graphs."   If $\gamma$ is unoriented at sufficiently many vertices, then $z(\gamma)$ lies in the stable range, and must be trivial in homology by Galatius's result.  In fact, as mentioned in the introduction, only two $z(\gamma)$ are known to be non-trivial in homology.  These arise from the graph with two vertices and three edges and the graph with two vertices and five edges (with orientations at both vertices); the status of the rest is unknown.  It is possible to extend the methods of this paper to give a proof that many graphs $\gamma$ 
with unoriented vertices give rise to trivial classes, including many that don't lie in the stable range. 
    Based on these remarks, it is reasonable to conjecture that $z(\gamma)$ is a nontrivial homology class if and only if all vertices of $\gamma$ are oriented.


\section{Relation to the Morita cocycle}
In this section we show that the cycles $z(\gamma)$ correspond to the cocycles $\mu_\gamma$ described in \cite{ConVog04}.  First we review the definition of these cocycles in the present context.

Let $\mathcal G$ denote the vector space spanned by admissible vertex-oriented graphs, modulo the relations $(G,or)+(G,-or)=0$, where $or$ denotes the total vertex orientation.  (A consequence of these relations is that a graph with an orientation-reversing automorphism is zero in $\mathcal G$.)
The \emph{Morita cocycle} is a linear map  $\mu:\mathfrak C^*\to \mathcal G$, defined as follows.

On a generator $(G, F)$, the Morita cocycle is zero unless $(G, F)$ has the following very special form:
\begin{itemize}
\item $G$ is trivalent.
\item $G$ decomposes as a subgraph $C$, which is the union of disjoint odd-sided simple polygons $C_i$, together with edges $E$  connecting these polygons. 
\item $F$ is a maximal forest in $C$.  
\end{itemize}
In this case $\mu(G,F)$ is equal to the graph $\gamma$ formed by collapsing each component $C_i$ of $C$ to a  vertex of $\gamma$ (see Figure 2). The orientation at the vertex of $\mu(G,F)$ corresponding to $C_i$ is induced by the orientation of $(G,F)$.  Specifically,  we have $(G, F)=\pm (G, F')$, where $F'$ has the same edges as $F$, ordered so that the edges in each linear component $F\cap C_i$ are numbered consecutively along the line.  This induces a cyclic orientation of $C_i$, which can be used to orient the incoming edges of the vertex in $\mu(G,F')=\pm \mu(G,F)$.  

\begin{center}
\includegraphics[width=4in]{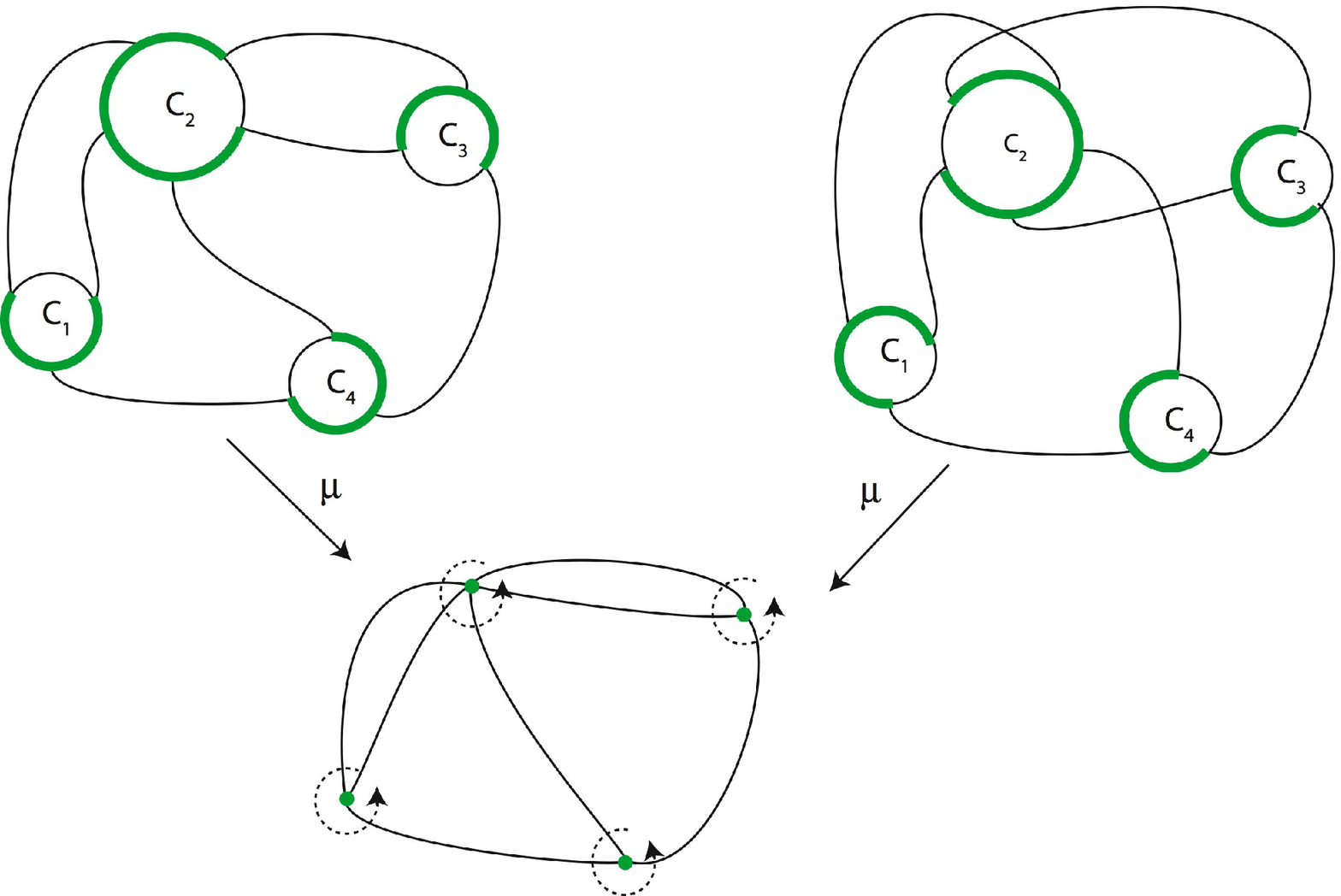}
\centerline{Figure 2:  Two graphs (G,F) with the same image under $\mu$}
\end{center}

The Morita cocycle takes values in a large vector space, $\mathcal G$. Projecting $\mathcal G$ onto $1$ dimensional subspaces produces infinitely many cocycles as follows.

\begin{definition}Let $0\neq{\bf g}\in \mathcal G$. The cocycle  $\mu_{\bf g}\colon \mathfrak C^*\to \Q$ is defined to be $\mu$ composed with the projection onto the subspace generated by ${\bf g}$.  In particular, we get a cocycle $\mu_\gamma$ with values in $\Q$ for every admissible vertex-oriented graph $\gamma$ with no orientation-reversing automorphisms.
 \end{definition}

To state the following theorem we need some preliminaries.  Let $\{,\}$ denote the canonical pairing $V^*\otimes V\to \mathbb Q$ of a vector space  with its dual, i.e. $\{ \phi,\vec v\}= \phi(\vec v)$. We have chain isomorphisms
$$\spine{p}{q}\to\left(\spine{p}{q}\right)^{**}\to\left(\dspine{p}{q}\right)^*$$
The left isomorphism is the canonical one between a finite dimensional vector space and its double dual, sending $\vec v$ to the functional $\psi_{\vec v}$ determined by  $\{\psi_{\vec v},\phi\}=\phi(\vec v)$.
The right isomorphism is induced by the isomorphism
$$\dspine{p}{q}\to\left(\spine{p}{q}\right)^*$$ given by the chain map
$X\mapsto \langle \cdot, X\rangle$.  We denote the composition by $\eta\colon \spine{p}{q}\to\left(\dspine{p}{q}\right)^*$.

\begin{theorem}
For every admissible vertex-oriented graph $\gamma$, we have $$c_\gamma \cdot \mu_\gamma= \eta(z(\gamma)),$$ where $c_\gamma$ is a nonzero integer.
\end{theorem}
\proof{

Analyzing what happens to $z(\gamma)$ in the above diagram, we see that
$z(\gamma)$ is mapped to a functional $f_\gamma$ defined by
$\{f_\gamma,Y\}=\langle z(\gamma),Y\rangle$.
The theorem will follow once we show that
$$\langle z(\gamma),X\rangle=c_\gamma\cdot\mu_\gamma(X)$$
where $$c_\gamma=|\operatorname{Aut}(\gamma)|\prod_{v}(2|v|) $$
the product ranging over all  vertices of $\gamma$ and $|v|$ denoting the valence of $v$.
To see this, first note that $\operatorname{Aut}(\gamma)$ acts by permuting the summands of $z(\gamma)$. Secondly, if two summands are isomorphic, it is clear that the isomorphism must correspond to an automorphism of $\gamma$. Thirdly, If this corresponding automorphism of $\gamma$ is the identity, this means that the permutations at each vertex actually differ by elements of the dihedral group. Thus we can conclude each term $X$ in $z(\gamma)$ appears $\displaystyle
\frac{|\operatorname{Aut}(\gamma)|}{|\operatorname{Aut}(X)|} \prod_{v}(2|v|)$ times, from which the desired conclusion follows.}

In \cite{ConVog04} we defined the Morita cocycle $\mu^{\text{old}}_\gamma$ on the {\it forested graph complex}, $f\mathcal{G}$.  The superscript refers to the fact that this is the ``old" definition.
This forested graph complex is obtained by considering the vertical filtration of the double complex $\mathfrak C^*$.  
The columns are exact except at the top (where the graphs $G$ are trivalent), so the $E^1$ page of the associated spectral sequence collapses to a single row, with terms $\mathfrak C^{p,2n-2}/\operatorname{im}(\delta_C)$ and differentials induced by $\delta_R$.  This is exactly the forested graph complex, with $\operatorname{IHX}$-relations exactly corresponding to $\operatorname{im}(\delta_C)$. Thus, we have a quasi-isomorphism $\xi\colon\mathfrak C^*\to f\mathcal{G}$ defined by passing to the second page of the spectral sequence.

The next proposition states that the Morita cocycle we defined in \cite{ConVog04} is the same as the cocycle $\mu$ defined in this section.

\begin{proposition}  The Morita cocycle, $\mu^{\text{old}}_\gamma$, is equal to the functional $\mu_\gamma\circ \xi$. 
\end{proposition}

\proof{This follows straight from the definitions.}



\section{Stable Triviality}

If we choose a basepoint $b$ on $C$, the cycle $z(\gamma)$ defined in Section 3 becomes a cycle $z(\gamma)$ in the quotient of the spine $L_n$ of Auter space.  In Figure 3 we indicate $z(\gamma)$ in a neighborhood of $C_1$, where the shading is meant to represent the sum over all maximal forests in the shaded circle.

\begin{center}
\includegraphics[width=2.2in]{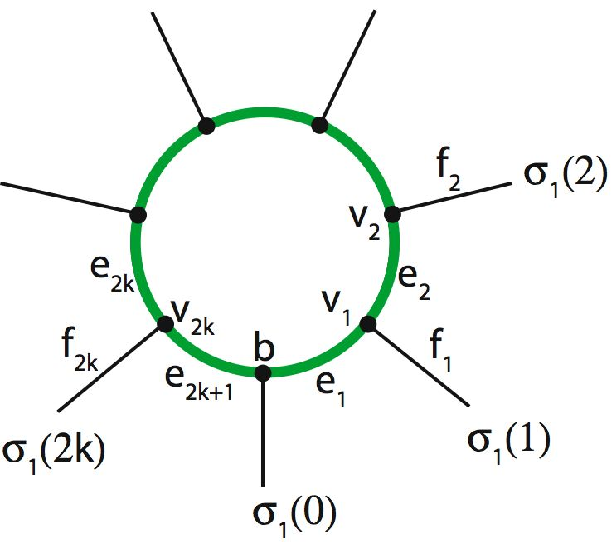}
\centerline{Figure 3: $z(\gamma)$ near $C_1$}\label{Zgamma}
\end{center}

The vertices and edges have been given names, for convenience, where the edges of $C_1$ have been numbered consistently with the definition of $z(\gamma)$ coming from Proposition \ref{dR}. Finally, the fact that the $f_j$ edges attach to the rest of the graph via the permutation $\sigma_1$ has been indicated.

 The stabilization map $s\colon H_i(\Aut{n})\to H_i(\Aut{n+1})$ is induced by the map $L_n\to L_{n+1}$ which adds a loop to the basepoint of each marked graph $(G,g)$.
In this section we show that the image of $z(\gamma)$ under this map is a boundary in the quotient of $L_{n+1}$, i.e. the homology stabilization map  sends $z(\gamma)$ to zero.  

We first establish some notation. If $G$ is any graph,   $b$ is the basepoint of $G$  
and $e$ is any edge of $G$,  let $G[e]$ denote the graph obtained from $G$ by adding an edge from $b$ to the midpoint of $e$.  This divides the edge $e$ into two new edges and increases the rank of $G$ by one.  Similarly, for any vertex $v$, let $G[v]$ denote the graph obtained from $G$ by adding an edge from $b$ to $v$.
Finally,  if $A$ is a subgraph of $G$, let    $\bdry A$  denote the boundary of $A$, i.e. the set of edges of $G$ with one vertex in $A$ and one vertex in $G-A$.

\begin{theorem} For any admissible vertex-oriented graph $\gamma$, and any choice of basepoint, the image of $z(\gamma)$ under the stabilization map $$ H_i(\Aut{n})\to~H_i(\Aut{n+1})$$
   is null-homologous.
\end{theorem}

\proof{Recall that the edges of each $C_i$ are ordered cyclically, and the basepoint of $G$ is on $C_1$, between the first and last edges of that polygon.  Let $Z^+=s(z(\gamma))$ be the stabilization of $z(\gamma)$, i.e. $Z^+$ is equal to  $z(\gamma)$ with loops added at the basepoints of all graphs.  We need to produce a chain $W$ with $\bdry(W)=Z^+$.

The polygon $C_1$ has $n=2k+1$ edges  $e_1,\ldots,e_{n}$.  If we divide $e_i$ into two edges, we obtain a new polygon $C_1^{(i)}$ with $n+1=2k+2$ edges, ordered cyclically so that the two halves of $e_i$ become $i$-th and ${(i+1)}$-st edges in $C_1^{(i)}$.   Set $C^{(i)}=C_1^{(i)}\cdot C_2\cdot\ldots\cdot C_k$.

Define  $$X_i=\sum_\sigma\sum_{F' } \epsilon_\sigma\epsilon_{F'}(G^\sigma[e_i], F'),$$
where $F'$ runs over all maximal subforests of $C^{(i)}$. Since all of the ``action" in what follows will be at $C_1^{(i)}$, we represent $X_i$ in Figure 4 by showing $C_1^{(i)}$ in green and indicating $\sigma_1$ at the tips of the edges emanating from $C_1^{(i)}$.

\begin{center}
\includegraphics[width=2in]{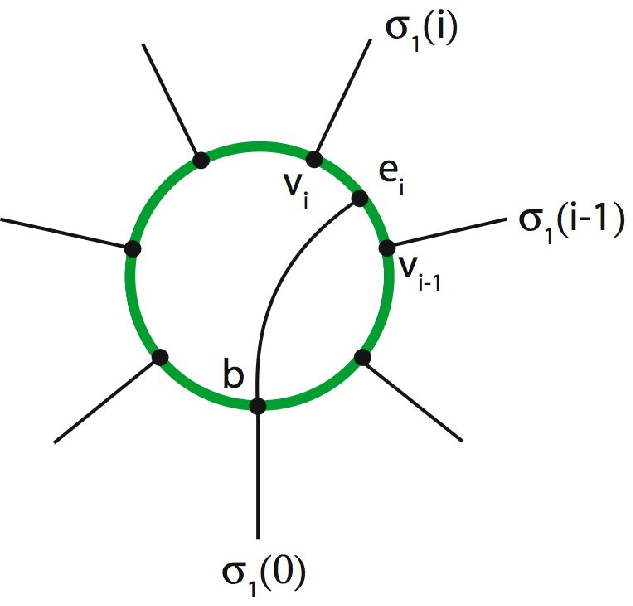}
\centerline{Figure 4: $X_i$ near $C_1$}
\end{center}

We have $\displaystyle d_R(X_i)=\sum_\sigma\epsilon_\sigma\,d_R(\sum_{F' } \epsilon_{F' }(G^\sigma[e_i], F'))=0$ by Proposition~\ref{dR}.

\medskip

Let $v_0=b, v_1,\ldots, v_{2k}$ be the  vertices of $C_1$, and $f_i$ the edge in $\bdry C_1$ attached at $v_i$. Let $m_i$ denote the midpoint of $f_i$ and define
$$Y_i=\sum_\sigma \sum_F\epsilon_\sigma \epsilon_F (G^\sigma[f_i], F\cdot [v_i,m_i])$$
where  $F$ runs through all maximal forests in $C$
and $$Y'_i=\sum_\sigma \sum_F\epsilon_\sigma \epsilon_F(G^\sigma[f_i], F\cdot [b,m_i]).$$
(see Figure 5)
\begin{center}
\includegraphics[width=5in]{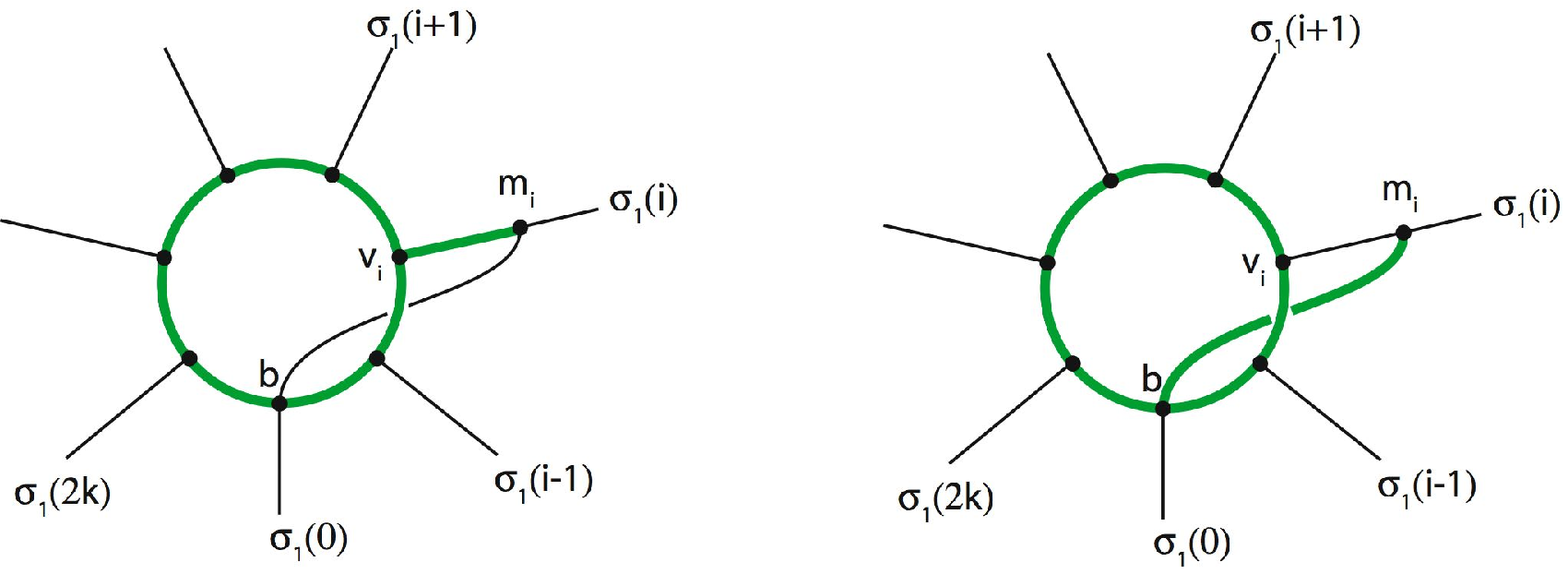}
\centerline{Figure 5: $Y_i$ and $Y_i'$ near $C_1$}\label{Yi}
\end{center}

Then
\begin{align*}
d_R(Y_i)&=\sum_\sigma \sum_F\epsilon_\sigma \epsilon_F(G^\sigma[f_i],F)+\sum_\sigma \sum_F\epsilon_\sigma \epsilon_F(G^\sigma[f_i],d_R(F)\cdot [b,m_i])\\
&=\sum_\sigma\sum_F \epsilon_\sigma\epsilon_F(G^\sigma[f_i],F),
\end{align*}
 where the second summation term is zero by Proposition \ref{dR} and, similarly, 
\begin{align*}
d_R(Y'_i)&=\sum_\sigma \sum_F\epsilon_\sigma \epsilon_F(G^\sigma[f_i],F)+\sum_\sigma \sum_F\epsilon_\sigma \epsilon_F(G^\sigma[f_i],d_R(F)\cdot [v_i,m_i])\\
&=\sum_\sigma\sum_F \epsilon_\sigma\epsilon_F(G^\sigma[f_i],F)
\end{align*}

Thus $d_R(Y_i-Y'_i)=0$.

Next we compute $d_C$ of $X_i, Y_i$ and $Y'_i$. To simplify notation, we define the following three chains 
(see Figure 6). 

\begin{center}
\includegraphics[width=6in]{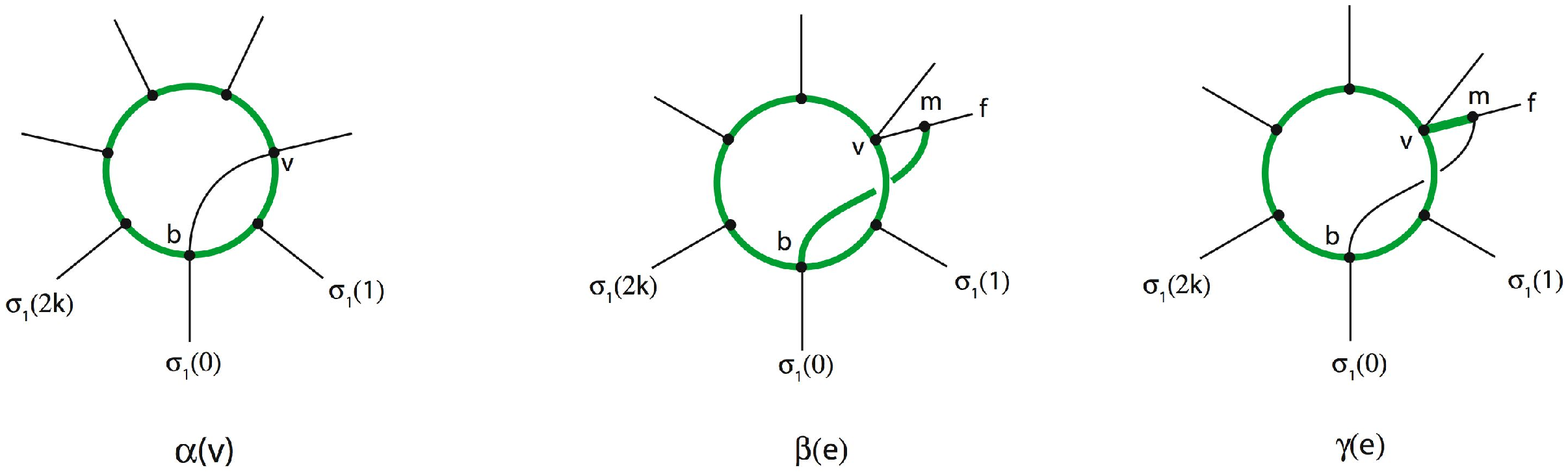}
\centerline{Figure 6: The chains $\alpha(v)$, $\beta(e)$ and $\gamma(e)$}\label{abc}
\end{center}

For any vertex $v$ of $C_1$, we set $$\alpha(v)=  \sum_\sigma \sum_F \epsilon_\sigma \epsilon_F (G^\sigma[v],F),$$ where $F$ runs over all maximal subforests of $C$.  In particular, $\alpha(b)=Z^+$.  For any edge $e$ of $C_1$, set 
$$\beta(e)=\sum_\sigma\sum_{F_e} \epsilon_\sigma\epsilon_{F_e} (G_e^{\sigma}[f], F_e\cdot [b,m]),$$
where $f$ is the boundary edge attached at the end of $e$, $m$ is the midpoint of $f$, $[b,m]$ is the edge from $b$ to $m$, and $F_e$ runs over all maximal forests in the image $C_e$ of $C$ in $G_e$ and
$$\gamma(e)=\sum_\sigma\sum_{F_e} \epsilon_\sigma\epsilon_{F_e} (G_e^{\sigma}[f], F_e\cdot [v,m]),$$ where $f$ is the boundary edge attached at the end of $e$, $m$ is the midpoint of $f$, $[v,m]$ is the edge from the endpoint of $e$ to $m$, and $F_e$ runs over all maximal forests in the image $C_e$ of $C$ in $G_e$.  

\begin{lemma}  
 $d_C(X_i)= (-1)^i [\alpha(v_{i})-\alpha(v_{i-1})].$  In particular, $d_C(X_1)=-\alpha(v_1) + Z^+.$ \end{lemma}
 \proof{ 
All terms other than $\alpha(v_{i})$ and $\alpha(v_{i-1})$ cancel upon symmetrization. 
(That is, as in the proof of Proposition \ref{dc}, there are two edges emanating from the same vertex of $G^\sigma_e$. When summing over all $\sigma$, this term will cancel with the term where these edges are transposed.)
See Figure 7.  
\begin{center}
\includegraphics[width=3.5in]{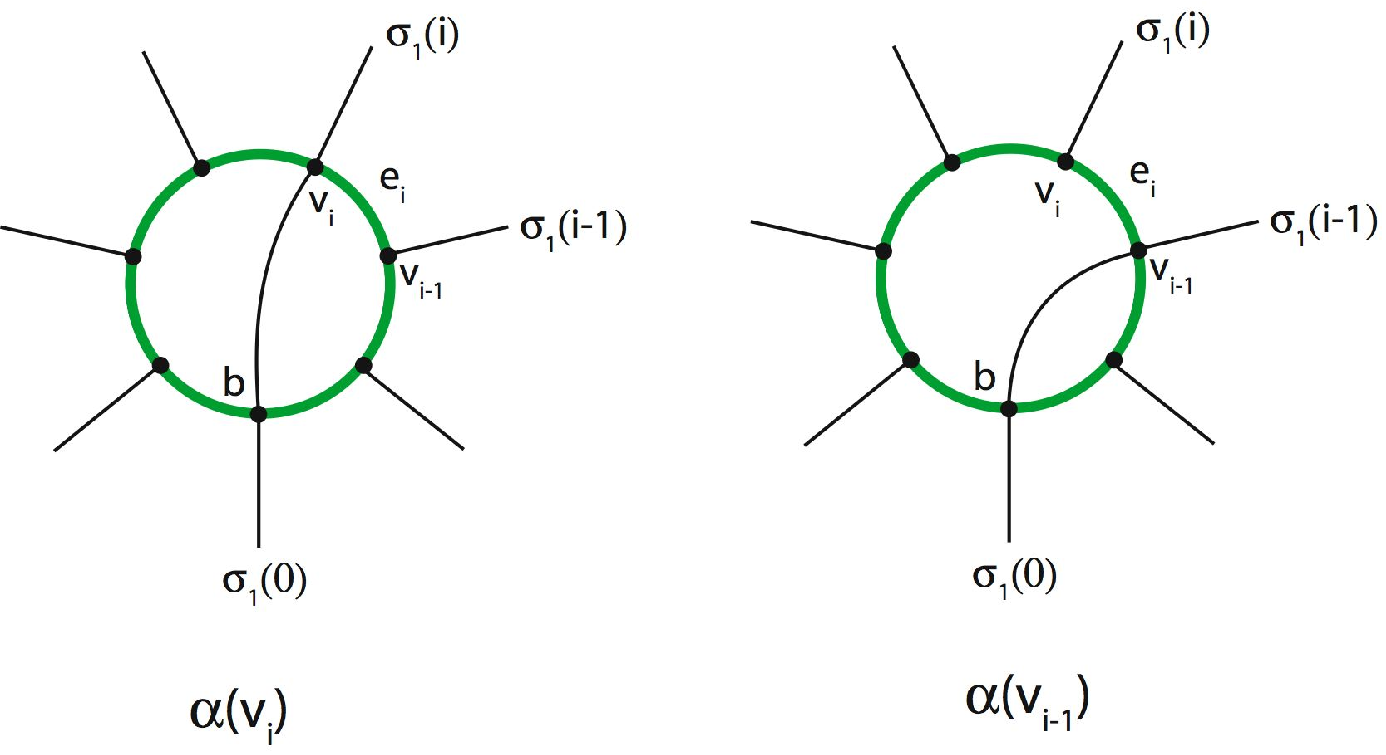}
\centerline{Figure 7: $d_C(X_i)$ near $C_1$}\label{dCXi}
\end{center}

 }

\begin{lemma} $ d_C(Y_i)= \alpha(v_i) + (-1)^i[ \gamma(e_{i+1})+\gamma(e_{i})] $.
 \end{lemma}
 \proof{  All terms other than $\alpha(v_i), \gamma(e_{i+1})$ and $\gamma(e_{i})$ cancel upon symmetrization.  See Figure 8.
\begin{center}
\includegraphics[width=6in]{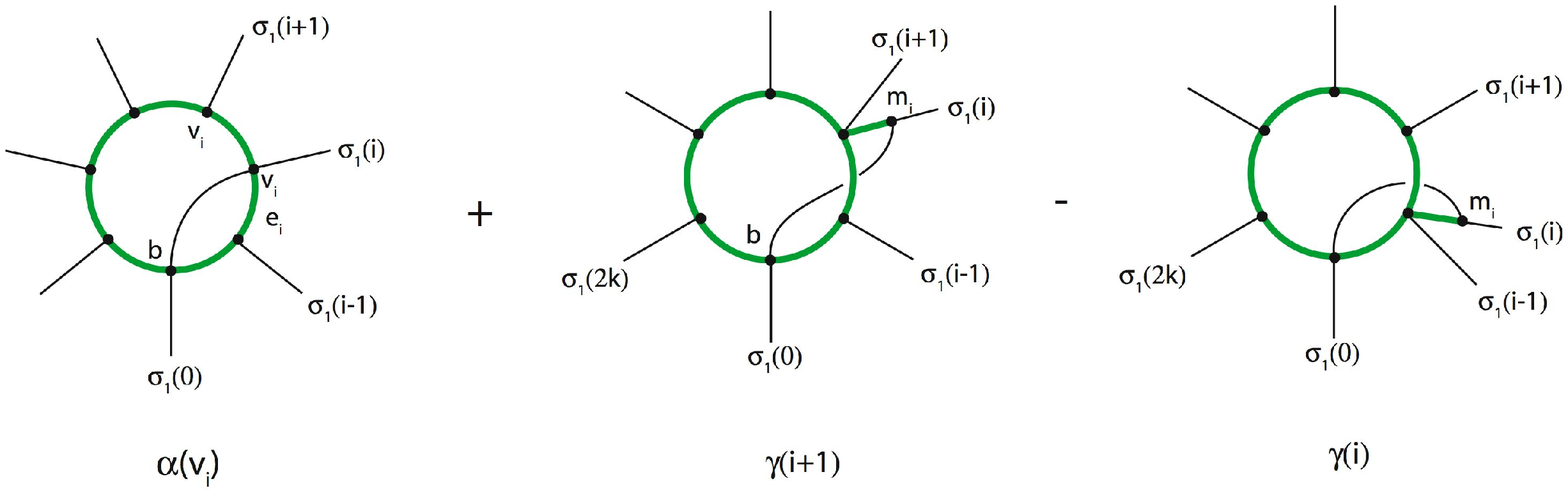}
\centerline{Figure 8: $d_C(Y_i)$ near $C_1$}\label{dCYi}
\end{center}

 }

\begin{lemma}$d_C(Y'_i)= (-1)^i[\beta(e_{i+1})+\beta(e_{i})]$.
\end{lemma}
\proof{  All terms other than $\beta(e_{i+1})$ and $\beta(e_{i})$ cancel upon symmetrization. See Figure 9. 
\begin{center}
\includegraphics[width=4in]{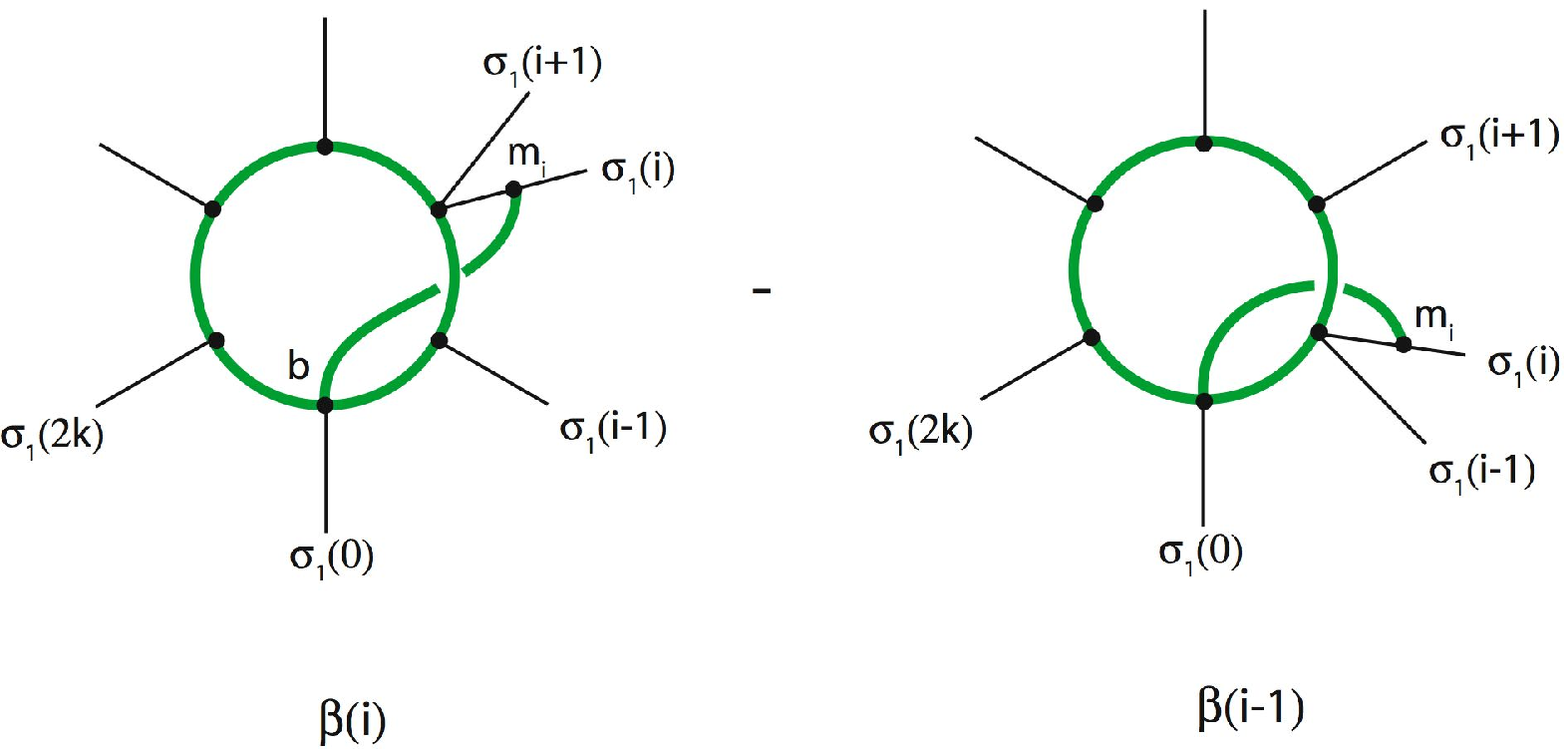}
\centerline{Figure 9: $d_C(Y_i')$ near $C_1$}\label{dCYiprime}
\end{center}

}

Now set $\displaystyle Y=\sum_{i=1}^{2k} Y_i$ and $\displaystyle Y'=\sum_{i=1}^{2k} Y_i'$.  All terms of $d_C(Y')$ cancel and we get $d_C(Y')=0$.  The terms of  $d_C(Y)$  also cancel except for the terms $\alpha(v_i)$, and we obtain  $\displaystyle d_C(Y)=\sum_{i=1}^{2k}\alpha(v_i)$. We remark that there is symmetry in the picture, and in fact $\alpha(v_i)=\alpha(v_{n-i})$, so we can write $\displaystyle d_C(Y)=2\sum_{i=1}^k \alpha(v_i)$.   

Finally, define
$$W=\frac{1}{2k} (Y-Y') + X_1 -\frac{k-1}{k}X_2 + \frac{k-2}{k}X_3 + \ldots + (-1)^{k+1} \frac{1}{k}X_k.$$  Our computations give $d_R(W)=0$ and 
\begin{align*}
d_C(W)&= \frac{1}{k}[\alpha(v_1)+\ldots+\alpha(v_k)]\\
& \ \ \ + [Z^+-\alpha(v_1)] + \frac{k-1}{k}[\alpha(v_1)-\alpha(v_2)] + \ldots + \frac{1}{k}[\alpha(v_{k-1})-\alpha(v_k)]\\
&=Z^+
\end{align*}
}
\def\cprime{$'$} \def\cprime{$'$}


\end{document}